\documentclass[12pt]{article}

\usepackage{hyperref}
\usepackage{xcolor}
\usepackage{amssymb}
\usepackage{amsfonts}
\usepackage{amsthm}
\usepackage{amsmath}
\usepackage{float}
\usepackage{bigints}
\usepackage{fullpage}
\usepackage{mathtools}
\usepackage[utf8]{inputenc} 
\usepackage{cancel}
\usepackage{verbatim}

\newtheorem{theorem}{Theorem}

\newtheorem{lemma}[theorem]{Lemma}

\newtheorem{conjecture}[theorem]{Conjecture}
\newtheorem{definition}{Definition}

\newtheorem{corollary}[theorem]{Corollary}
\newtheorem{remark}[theorem]{Remark}

\newcommand{\p}{\mathbf{p}}
\newcommand{\q}{\mathbf{q}}

\newcommand{\oo}{\mathbf{o}}

\newcommand{\B}{\mathbf{B}}

\newcommand{\Ee}{\mathbb{E}}

\newcommand{\Rr}{\mathbb{R}}

\newcommand{\Ed}{\Ee^d}

\newcommand{\iprod}[2]{\left<#1,#2\right>}

\newcommand{\noshow}[1]{}
\newcommand{\ivol}[2][k]{{\rm V}_{#1}\left(#2\right)}

\title{Notes on uniform contractions of congruent balls \footnote{\hspace*{.35cm} 2010 Mathematics Subject Classification: 52A20, 52A22.}}

\author{K\'{a}roly Bezdek\thanks{Partially supported by a Natural Sciences and 
Engineering Research Council of Canada Discovery Grant.}
}

\begin{document}

\maketitle

\begin{abstract}
In this note, we obtain stronger results on uniform contractions of intersections and unions of congruent balls, thereby extending previous partial confirmations of the Kneser–Poulsen conjecture in these settings. In addition, we establish Blaschke--Santal\'o-type inequalities for $r$-ball bodies.
\end{abstract}

\section{Introduction}\label{sec:intro}

Let $\Ee^d$, $d\geq 1$ denote the $d$-dimensional Euclidean space with $\oo\in \Ee^d$ standing for the origin. We denote the Euclidean norm of a vector $\p$ in $\Ee^d$ by $|\p|:=\sqrt{\iprod{\p}{\p}}$, where $\iprod{\cdot}{\cdot}$ is the 
standard inner product. The closed Euclidean ball of radius $r\geq 0$ centered at $\p\in\Ed$ is denoted by $\B^d[\p,r]:=\{\q\in\Ed\ |\  |\p-\q|\leq r\}$. It will be convenient to use for sets $A, B\subset \Ed$ the notations $A+B:=\{\mathbf{a}+\mathbf{b}\ |\ \mathbf{a}\in A\ {\rm and}\ \mathbf{b}\in B\}$ and $\lambda A:=\{\lambda\mathbf{a}\ |\  \mathbf{a}\in A\}$ for given $\lambda\in \Rr$. Let $A\subset\Ed$ be a compact convex set, and 
$0\leq k\leq d$. We denote the {\it $k$-th intrinsic volume} of $A$ by $\ivol{A}$. It is well known that $\ivol[d]{A}$ is the $d$-dimensional 
volume (Lebesgue measure) of $A$, $2\ivol[d-1]{A}$ is the surface area of $A$, $\frac{2\omega_{d-1}}{d\omega_d}\ivol[1]{A}$ is equal to the mean width of $A$, and $\ivol[0]{A}:=1$ is the Euler characteristic of $A$, where $\omega_d$ stands for the volume of a $d$-dimensional unit ball, that is, $\omega_d=\frac{\pi^{\frac{d}{2}}}{\Gamma(1+\frac{d}{2})}$. For properties of intrinsic volumes, including Steiner's formula
\begin{equation}\label{Steiner}
\ivol[d]{A+\epsilon\B^d[\oo,1]}=\sum_{i=0}^{d}\omega_{d-i}\ivol[i]{A}\epsilon^{d-i},
\end{equation}
we refer the interested reader to \cite{GHSch}. In this note, for simplicity $\ivol{\emptyset}=0$ for all $0\leq k\leq d$.

Following \cite{BeNa}, we say that a (labeled) point set $Q:=\{\q_1, \dots , \q_N\}\subset\Ee^d$ is a {\it uniform contraction} of a (labeled) point set $P:=\{\p_1, \dots ,\p_N\}\subset\Ee^d$ with {\it separating value} $\lambda>0$ in $\Ee^d$  if
\begin{equation}\label{uniform}
|\q_i-\q_j|\leq\lambda\leq |\p_i-\p_j|\ {\rm holds}\ {\rm for}\  {\rm all}\ 1\leq i<j\leq N.
\end{equation} 
Furthermore,  the (labeled) point set $Q:=\{\q_1, \dots , \q_N\}\subset\Ee^d$ is a {\it uniform contraction} of the (labeled) point set $P:=\{\p_1, \dots ,\p_N\}\subset\Ee^d$ in $\Ee^d$ if there exists a $\lambda>0$ such that $Q$ is a uniform contraction of $P$ with separating value $\lambda$ in $\Ee^d$. For a compact set $\emptyset\neq A\subset \Ee^d$ let the {\it circumradius} ${\rm cr}(A)$ of $A$ be defined by ${\rm cr}(A):= \inf \{ r>0\ |\ A\subseteq \B^d[\p,r]\ {\rm for}\ \p\in \Ee^d\}$. If $\emptyset\neq A\subset \Ee^d$ is a compact set with $0\leq {\rm cr}(A)\leq r$ for some $r>0$,
let $A^r:=\bigcap\{\B^d[\p,r]\ |\ \p\in A\}$ be the  the {\it $r$-dual of $A$} (\cite{Be18}) or equivalently, the {\it $r$-ball body generated by $A$} (\cite{Be22}, \cite{Be23}), which for $r=1$ is simply called a {\it ball-body generated by $A$}. If $A$ is a finite set, then we call $A^r$ an {\it $r$-ball polyhedron} or simply, a {\it ball-polyhedron} (\cite{BLNP}) for $r=1$.
Now we can recall Theorem 1.4 of  \cite{BeNa} as follows: Let $r>0$, $d>1$, $1\leq k\leq d$, and let $N\geq (1+\sqrt{2})^d$. If $Q:=\{\q_1, \dots , \q_N\}\subset\Ee^d$ is a uniform contraction of $P:=\{\p_1, \dots ,\p_N\}\subset\Ee^d$ in $\Ee^d$, then $\ivol[k]{P^r}\leq \ivol[k]{Q^r}$. In Section~\ref{first proof} of this note we give a proof of the following stronger result.

\begin{theorem}\label{KP-type improved}
Let $r>0$, $d>1$, $1\leq k\leq d$, and let $N\geq \left(1+\sqrt{\frac{2d}{d+1}}\right)^d$. If $Q:=\{\q_1, \dots , \q_N\}\subset\Ee^d$ is a uniform contraction of $P:=\{\p_1, \dots ,\p_N\}\subset\Ee^d$ in $\Ee^d$, then 
\begin{equation}
\ivol[k]{P^r}\leq \ivol[k]{Q^r}.
\end{equation} 
\end{theorem}

Recall that a (labeled) point set $Q:=\{\q_1, \dots , \q_N\}\subset\Ee^d$ is a {\it contraction} of a (labeled) point set $P:=\{\p_1, \dots ,\p_N\}\subset\Ee^d$ in $\Ee^d$  if
$|\q_i-\q_j|\leq |\p_i-\p_j|\ {\rm holds}\ {\rm for}\  {\rm all}\ 1\leq i<j\leq N$. In 1954, E. T. Poulsen \cite{Po54} and in 1955, M. Kneser \cite{Kn55} independently conjectured the following inequalities for the case when $r_1=\dots=r_N$. If $Q:=\{\q_1, \dots , \q_N\}\subset\Ee^d$ is a contraction of $P:=\{\p_1, \dots ,\p_N\}\subset\Ee^d$ in $\Ee^d$, then according to the general form of the Kneser--Poulsen conjecture the following inequalities hold for all $r_1>0, r_2> 0, \dots, r_N> 0$:
\begin{equation}\label{K-P}
\ivol[d]{\bigcup_{i=1}^{N}\mathbf{B}^d[{\bf p}_i, r_i]}\ge \ivol[d]{\bigcup_{i=1}^{N}\mathbf{B}^d[{\bf q}_i, r_i]}
\ {\rm and}\ 
\ivol[d]{\bigcap_{i=1}^{N}\mathbf{B}^d[{\bf p}_i, r_i]}\le \ivol[d]{\bigcap_{i=1}^{N}\mathbf{B}^d[{\bf q}_i, r_i]}.
\end{equation}
In \cite{BeCo}, the author and Connelly proved the Kneser--Poulsen conjecture (i.e., (\ref{K-P})) in $\Ee^2$. The conjecture remains open in dimension 3 and higher. For more details on the Kneser--Poulsen conjecture we refer the interested reader to the recent survey papers \cite{BLN} and \cite{Cs18}. 

It seems that Alexander~\cite{Alex} was the first who conjectured that under an arbitrary contraction of the center points of finitely many congruent disks in the plane, the perimeter of the intersection of the disks cannot decrease. In other words, if $Q:=\{\q_1, \dots , \q_N\}\subset\Ee^2$ is a contraction of $P:=\{\p_1, \dots ,\p_N\}\subset\Ee^2$ for $N>1$ in $\Ee^2$, then $\ivol[1]{P^r}\leq \ivol[1]{Q^r}$ holds for any $r>0$. If true, it would give a sharpening of the Kneser-Poulsen conjecture for the special case of the intersection of congruent disks. Alexander's conjecture was proved for $N=2,3,4$ in \cite{BeCoCs}. As $\left(1+\sqrt{\frac{4}{3}}\right)^2=4.6427...$ therefore Theorem ~\ref{KP-type improved} with $d=2$ and $k=1$ proves Alexander's conjecture for any uniform contraction with $N\geq 5$. Thus, Alexander's conjecture is proved for all uniform contractions:

\begin{corollary}\label{Alex-type in 2D}
Let $Q:=\{\q_1, \dots , \q_N\}\subset\Ee^2$ be a uniform contraction of $P:=\{\p_1, \dots ,\p_N\}\subset\Ee^2$ for $N>1$ in $\Ee^2$. Then $\ivol[1]{P^r}\leq \ivol[1]{Q^r}$ holds for any $r>0$.
\end{corollary}

Recall that the author and Connelly (Corollary 4 in ~\cite{BeCo}) proved the Kneser--Poulsen conjecture (i.e., (\ref{K-P})) for all $1\leq N\leq d+3$ in $\Ee^d, d>2$. Thus Theorem~\ref{KP-type improved} for $k=d$ can be strengthened as follows.

\begin{corollary}\label{KP-style strengthened}
Let $r>0$, $d>2$, and let either $1\leq N\leq d+3$ or $N\geq \left(1+\sqrt{\frac{2d}{d+1}}\right)^d$. If $Q:=\{\q_1, \dots , \q_N\}\subset\Ee^d$ is a uniform contraction of $P:=\{\p_1, \dots ,\p_N\}\subset\Ee^d$ in $\Ee^d$, then 
\begin{equation}\label{KP-style strengthened formula}
\ivol[d]{P^r}\leq \ivol[d]{Q^r}.
\end{equation}
In particular, (\ref{KP-style strengthened formula}) holds for all uniform contractions of balls of radius $r$ in $\Ee^3$ if and only if it holds for all uniform contractions of $7, 8, 9, 10$, and $11$ balls of radius $r$ in $\Ee^3$.
\end{corollary}

If  $\emptyset\neq A\subset \Ee^d$ is a compact set with $0\leq {\rm cr}(A)\leq r$ for some $r>0$, then let the {\it $r$-convex hull}  ${\rm conv}_r(A)$ (also called $r$-ball convex hull (\cite{Be22}, \cite{Be23}) or simply, spindle convex hull (\cite{BLNP}) for $r=1$) be defined by ${\rm conv}_r(A):=\bigcap\{\B^d[\p,r]\ |\ A\subseteq \B^d[\p,r] \ {\rm for}\ \p\in \Ee^d\}$. It is easy to see that $A$ is an $r$-ball body if and only if $A={\rm conv}_r(A)$. Furthermore, Lemma 5 of \cite{Be18} implies that 
\begin{equation}\label{basic}
A^r=\left({\rm conv}_r(A)\right)^r.
\end{equation}
It was observed in \cite{Ca} as well as in \cite{BeCoCs} that Alexander's conjecture \cite{Alex} is equivalent to the following conjecture: If $P$ is a set of finitely many points lying in a circular disk of radius $r>0$ in $\Ee^2$, then the perimeter (resp., first intrinsic volume) of the $r$-convex hull of $P$ under arbitrary contraction of $P$ cannot increase in $\Ee^2$. Hence, it is natural to put forward the following extended version of Alexander's conjecture, which is supported also by Theorem~\ref{KP-type improved}.

\begin{conjecture}\label{Alexander extended}
Let $d>1$ and $r>0$. Under arbitrary contraction of the center points of finitely many balls of radius $r$ in $\Ed$, the mean width (resp., first intrinsic volume) of the intersection of the balls cannot decrease. Equivalently, if $P$ is a set of finitely many points lying in a ball of radius $r$ in $\Ed$, then the mean width (resp., first intrinsic volume) of the $r$-convex hull of $P$ under arbitrary contraction of $P$ cannot increase in $\Ed$.
\end{conjecture}

We note that the equivalence part of Conjecture~\ref{Alexander extended} follows from Lemma~\ref{difference identity} of Section~\ref{first proof} in a straightforward way. 

\begin{remark}\label{contractions of molecules}
Here we recall that an analogue of Theorem~\ref{KP-type improved} for unions of congruent balls have been proved in \cite{BeCo}, \cite{Be20} and can be stated as follows. Let $r>0$, $d>2$, and let either $1\leq N\leq d+3$ or $N\geq 2^d$. For a set $\emptyset\neq A\subset \Ee^d$ let  $A_r:=\bigcup\{\B^d[\p,r]\ |\ \p\in A\}$. If $Q:=\{\q_1, \dots , \q_N\}\subset\Ee^d$ is a uniform contraction of $P:=\{\p_1, \dots ,\p_N\}\subset\Ee^d$ in $\Ee^d$, then 
\begin{equation}\label{KP-style inequality for unions}
\ivol[d]{Q_r}\leq \ivol[d]{P_r}.
\end{equation}
In particular, (\ref{KP-style inequality for unions}) holds for all uniform contractions of balls of radius $r$ in $\Ee^3$ if and only if it holds for all uniform contractions of $7$ balls of radius $r$ in $\Ee^3$. 
\end{remark}

We strengthen the results of Remark~\ref{contractions of molecules} as follows. Let $\mathcal{P}_n:=\{\mathbf{B}^d[\mathbf{c}_i, 1]\  |\ 1\leq i \leq n\ {\rm with}\ |\mathbf{c}_j-\mathbf{c}_k|\geq 2\ {\rm for \ all}\ 1\leq j<k\leq n\}$ be a packing of $n>1$ unit balls in $\Ee^d$, $d>2$. Let $C:=\{\mathbf{c}_1,\dots ,\mathbf{c}_n\}$ and $\mu\geq 1$. Moreover, let
$$\delta_d(n,\mu):=\max_{\mathcal{P}_n}\frac{\ivol[d]{C_1}}{\ivol[d]{C_{\mu}}}=\max_{\mathcal{P}_n}\frac{n\omega_d}{\ivol[d]{\cup_{i=1}^{n}\mathbf{B}^d[\mathbf{c}_i, \mu]}}\ {\rm and}\ \delta_d(\mu):=\sup_{n}\delta_d(n,\mu).$$ Furthermore, let $\delta_d$ denote the {\it packing density} of unit balls in $\Ee^d$, which is the supremum of
the upper densities of all unit ball packings in $\Ee^d$. Next, we define the function $f(d, \mu)$ for each positive integer $d\geq 7$ over the reals $\mu\geq 1$ in four steps. Clearly,
\begin{equation}\label{step 1 of defining f}
f(d, \mu):= 1\geq \delta_d(\mu)
\end{equation}
holds for all $1\leq \mu<\sqrt[d]{d}$. Next, recall that according to Theorem 2 of \cite{BeLa15} 
\begin{equation}\label{step 2 of defining f}
f(d, \mu):=\frac{d+2}{2}{\mu}^{-d}\geq \delta_d(\mu)
\end{equation}
holds for all $\sqrt[d]{d}\leq\mu<\sqrt{\frac{2d}{d+1}}$. In order to do the next step let us take a regular $d$-dimensional simplex of edge length $2$ in $\Ee^d$ and then draw a $d$-dimensional unit ball around each vertex of the simplex. Let $\sigma_d$ denote the ratio of the volume of the portion of the simplex covered by balls to the volume of the simplex. Then according to the well-known result of Rogers \cite{Ro58}
\begin{equation}\label{step 3 of defining f}
f(d, \mu):=\sigma_d \geq \delta_d(\mu)
\end{equation}
holds for all $\sqrt{\frac{2d}{d+1}}\leq \mu<2$. Finally, the Lemma of \cite{Bez02} yields that
\begin{equation}\label{step 4 of defining f}
f(d, \mu):=\delta_d \geq \delta_d(\mu)
\end{equation}
holds for all $\mu\geq 2$.

\begin{theorem}\label{improved contractions of molecules}
Let $d\geq 7$, $\mu\geq 1$, and $N\geq f(d,\mu)(1+\mu)^d$. If $Q:=\{\q_1, \dots , \q_N\}\subset\Ee^d$ is a uniform contraction of $P:=\{\p_1, \dots ,\p_N\}\subset\Ee^d$ with separating value $\lambda$ in $\Ee^d$ and $r>\mu\frac{\lambda}{2}$, then
\begin{equation}\label{KP-style improved inequality for unions}
\ivol[d]{Q_r}\leq \ivol[d]{P_r}.
\end{equation}
\end{theorem}

\begin{remark}
We note that as $f(d,1)(1+1)^d=2^d$ and (\ref{KP-style improved inequality for unions}) holds trivially for all $r\leq\frac{\lambda}{2}$, therefore Remark~\ref{contractions of molecules} follows from Theorem~\ref{improved contractions of molecules}. Moreover, as (\ref{step 2 of defining f}) for $\sqrt[d]{d}\leq\mu<\sqrt{\frac{2d}{d+1}}$ implies that $f(d,\mu)(1+\mu)^d=\frac{d+2}{2}\left(1+\frac{1}{\mu}\right)^d<2^d$ holds for all sufficiently large $d$, therefore Theorem~\ref{improved contractions of molecules} yields a stronger result for $r>\mu\frac{\lambda}{2}$ than Remark~\ref{contractions of molecules}. Similarly, as (\ref{step 3 of defining f}) (resp., (\ref{step 4 of defining f})) for $\mu=\sqrt{\frac{2d}{d+1}}$ (resp., $\mu=2$) implies that $f(d,\mu)(1+\mu)^d=\sigma_d\left(1+\sqrt{\frac{2d}{d+1}}\right)^d<2^d$ (resp., $f(d,\mu)(1+\mu)^d=\delta_d3^d<2^d$) holds for all sufficiently large $d$, therefore Theorem~\ref{improved contractions of molecules} yields stronger results for $r>\mu\frac{\lambda}{2}$ than Remark~\ref{contractions of molecules}. Here we have used the upper estimates $\sigma_d<\frac{d}{e}2^{-0.499d}$ (\cite{Ro58}) and $\delta_d<2^{-0.599d}$ for sufficiently large d (\cite{CoEl}).
\end{remark}

In the rest of the paper we prove Theorems~\ref{KP-type improved} and~\ref{improved contractions of molecules}.

\section{Proof of Theorem~\ref{KP-type improved}}\label{first proof}

\begin{definition}
Let $r>0$. Let $d\geq 2$ and $1\leq l\leq d$ be integers. If $\emptyset\neq A\subset \Ee^d$ is a compact set with $0\leq{\rm cr}(A)\leq r$, then let the uniquely determined $R^{\Ee^d}_{r, l}(A)\geq 0$, defined by $\ivol[l]{{\rm conv}_r(A)}=\ivol[l]{\B^d[\oo, R^{\Ee^d}_{r, l}(A)]}$, be called the {\rm $l$-th intrinsic radius of $A$ relative to $r$ in $\Ee^d$}.
\end{definition}

Now, we are ready to state the following Blaschke--Santal\'o-type inequalities, which generalize Theorem 1 of \cite{Be18}, Theorem 1 of \cite{Be22}, and Theorem 1.2 of \cite{Be23}. (See also the Theorem of \cite{GHSch} for a special spherical analogue and Theorem 1.1 of \cite{PaPi} for a close stochastic analogue).

\begin{theorem}\label{BS-type}
Let $r>0$. Let $d\geq 2$ and $1\leq k\leq l\leq d$ be integers. If $\emptyset\neq A\subset \Ee^d$ is a compact set with $0<{\rm cr}(A)\leq r$ and $B:=\B^d[\oo, R^{\Ee^d}_{r, l}(A)]$, then
\begin{equation}\label{intrinsic}
\ivol[k]{A^r}\leq \ivol[k]{B^r}
\end{equation}
with equality if and only if ${\rm conv}_r(A)$ is a ball of radius $R^{\Ee^d}_{r, l}(A)$.
\end{theorem}

\begin{proof}
We start by recalling Proposition 2.5 of \cite{BeNa}. (See also Lemma 8 of \cite{Be20} for an extended version.)

\begin{lemma}\label{difference identity}
Let $r>0$ and let $d\geq 2$  be an integer. If $\emptyset\neq A\subset \Ee^d$ is a compact set with $0<{\rm cr}(A)\leq r$, then
\begin{equation}\label{difference is a ball}
A^r-{\rm conv}_r(A)=\B^d[\oo,r].
\end{equation}
\end{lemma}

Clearly, the Brunn--Minkowski inequality for intrinsic volumes (see for example, Eq. (74) in \cite{Gar}) combined with (\ref{difference is a ball}) of Lemma~\ref{difference identity} implies the following inequality.

\begin{corollary}\label{BM-type}
Let $r>0$ and let $d\geq 2$  be an integer. If $\emptyset\neq A\subset \Ee^d$ is a compact set with $0<{\rm cr}(A)\leq r$, then
\begin{equation}
\ivol[k]{A^r}^{\frac{1}{k}}+\ivol[k]{{\rm conv}_r(A)}^{\frac{1}{k}}\leq \ivol[k]{\B^d[\oo,r]}^{\frac{1}{k}}.
\end{equation}
\end{corollary}

Thus, Corollary~\ref{BM-type} and the positive homogeneity (of degree $k$) of $k$-th intrinsic volume imply in a straightforward way that

\begin{equation}\label{BM-type rearranged}
\ivol[k]{A^r}^{\frac{1}{k}}\leq r\ivol[k]{\B^d[\oo,1]}^{\frac{1}{k}}-\ivol[k]{{\rm conv}_r(A)}^{\frac{1}{k}}.
\end{equation}

Next, according to a well-known theorem of Alexandrov (see Eq. (20) of \cite{BuZa}) applied to the convex body ${\rm conv}_r(A)$ (which is by assumption a compact convex set with non-empty interior in $\Ee^d$), one obtains the following inequality
for the corresponding intrinsic volumes:
\begin{equation}\label {Al-type}
\left(\frac{\ivol[l]{{\rm conv}_r(A)}}{\ivol[l]{\B^d[\oo,1]}}\right)^k\leq \left(\frac{\ivol[k]{{\rm conv}_r(A)}}{\ivol[k]{\B^d[\oo,1]}}\right)^l ,
\end{equation}
where equality is attained if and only if  ${\rm conv}_r(A)$ is a ball.

Now, by definition $\ivol[l]{{\rm conv}_r(A)}=\ivol[l]{\B^d[\oo, R^{\Ee^d}_{r, l}(A)]}$ and therefore - using positive homogeneity of intrinsic volumes - (\ref{Al-type}) turns out to be equivalent to the simpler inequality
\begin{equation}\label{simpler Al-type}
\ivol[k]{\B^d[\oo, R^{\Ee^d}_{r, l}(A)]}\leq \ivol[k]{{\rm conv}_r(A)},
\end{equation}
where equality is attained if and only if ${\rm conv}_r(A)$ is a ball of radius $R^{\Ee^d}_{r, l}(A)$.

Finally, (\ref{BM-type rearranged}) and (\ref{simpler Al-type}) together with the positive homogeneity (of degree $k$) of $k$-th intrinsic volume imply in a straightforward way that
$$
\ivol[k]{A^r}^{\frac{1}{k}}\leq r\ivol[k]{\B^d[\oo,1]}^{\frac{1}{k}}-\ivol[k]{\B^d[\oo, R^{\Ee^d}_{r, l}(A)]}^{\frac{1}{k}}=\ivol[k]{\B^d[\oo, r-R^{\Ee^d}_{r, l}(A)]}^{\frac{1}{k}}=\ivol[k]{B^r}^{\frac{1}{k}}.
$$

Here, we have equality if and only if equality is attained in (\ref{BM-type rearranged}) and (\ref{simpler Al-type}), i.e., ${\rm conv}_r(A)$ is a ball of radius $R^{\Ee^d}_{r, l}(A)$. This completes the proof of Theorem~\ref{BS-type}.
\end{proof}

\begin{remark}
It is worth restating Theorem~\ref{BS-type} for $r$-ball bodies. Let $r>0$. Let $d\geq 2$ and $1\leq k\leq l\leq d$ be integers. If $\emptyset\neq \mathbf{A}\subset \Ee^d$ is an $r$-ball body different from a point (i.e., $\emptyset\neq\mathbf{A}\subset \Ee^d$ is a compact set with $0<{\rm cr}(\mathbf{A})\leq r$ and $\mathbf{A}={\rm conv}_r(\mathbf{A})$) and $\ivol[l]{\mathbf{A}}=\ivol[l]{\mathbf{B}}$ for $\mathbf{B}:=\B^d[\oo, R^{\Ee^d}_{r, l}(A)]$, then $\ivol[k]{\mathbf{A}^r}\leq \ivol[k]{\mathbf{B}^r}$ with equality if and only if $\mathbf{A}$ is a ball. In other words, the $k$-th intrinsic volume of $r$-duals of $r$-ball bodies of given $l$-th intrinsic volume in $\Ee^d$ is maximal only for balls. In particular, the $k$-th intrinsic volume of $r$-duals of $r$-ball bodies of given $d$-dimensional volume in $\Ee^d$ is maximal only for balls (Theorem 1.2 in \cite{Be23}).
\end{remark}

Now, we are ready to proceed with the proof of Theorem~\ref{KP-type improved}. We note that if $r<{\rm cr}(P)$, then clearly $\ivol[k]{P^r}=\ivol[k]{\emptyset}=0\leq \ivol[k]{Q^r}$, finishing the proof of Theorem~\ref{KP-type improved} in this case. Hence, for the rest of our proof we may assume that
\begin{equation}\label{circumradius}
{\rm cr}(P)\leq r .
\end{equation}

It follows from (\ref{basic}) and (\ref{circumradius}) that 
\begin{equation}\label{new representation}
\emptyset\neq P^r=\left(P_{\frac{\lambda}{2}}\right)^{r+\frac{\lambda}{2}}=\left({\rm conv}_{r+\frac{\lambda}{2}}\left(P_{\frac{\lambda}{2}}\right)\right)^{r+\frac{\lambda}{2}},
\end{equation}
where $P_{\frac{\lambda}{2}}:=\bigcup_{i=1}^N \B^d\left[\p_i,\frac{\lambda}{2}\right]$. By assumption on $P$, we get that $\{\B^d\left[\p_i,\frac{\lambda}{2}\right] | 1\leq i\leq N\}$ is a packing of balls in $\Ee^d$ and therefore
\begin{equation}\label{packing property}
\ivol[d]{{\rm conv}_{r+\frac{\lambda}{2}}\left(P_{\frac{\lambda}{2}}\right)}\geq \ivol[d]{\B^d\left[\oo, N^{\frac{1}{d}}\frac{\lambda}{2}\right]}.
\end{equation}
Hence, (\ref{new representation}) and (\ref{packing property}) combined with Theorem~\ref{BS-type} (for $1\leq k\leq l=d$) yield via homogeneity and monotonicity of intrinsic volume that
\begin{equation}\label{part 1}
\ivol[k]{P^r}\leq \ivol[k]{\B^d\left[\oo, r-\left(N^{\frac{1}{d}}-1\right)\frac{\lambda}{2}\right]}=\left(r-\left(N^{\frac{1}{d}}-1\right)\frac{\lambda}{2}\right)^k \ivol[k]{\B^d\left[\oo, 1\right]},
\end{equation}
where $r-\left(N^{\frac{1}{d}}-1\right)\frac{\lambda}{2}\geq 0$. The latter inequality and the assumption $N\geq \left(1+\sqrt{\frac{2d}{d+1}}\right)^d$ imply
\begin{equation}\label{useful condition}
0\leq r-\left(N^{\frac{1}{d}}-1\right)\frac{\lambda}{2}\leq r- \sqrt{\frac{2d}{d+1}}\frac{\lambda}{2}.
\end{equation}

Next, notice that $Q$ is a set of diameter at most $\lambda$ in $\Ee^d$ and therefore using Jung's theorem (\cite{De}) we get that ${\rm cr}(Q)\leq \sqrt{\frac{2d}{d+1}}\frac{\lambda}{2}$. This inequality and (\ref{useful condition}) easily imply that 
$\B^d\left[\mathbf{x}, r- \sqrt{\frac{2d}{d+1}}\frac{\lambda}{2}\right]\subseteq Q^r$ holds for some $\mathbf{x}\in \Ee^d$ with $ r- \sqrt{\frac{2d}{d+1}}\frac{\lambda}{2}\geq 0$. Thus, homogeneity and monotonicity of intrinsic volume yield
\begin{equation}\label{part 2}
\left(r- \sqrt{\frac{2d}{d+1}}\frac{\lambda}{2}\right)^k \ivol[k]{\B^d[\mathbf{x}, 1]}=\ivol[k]{\B^d\left[\mathbf{x}, r- \sqrt{\frac{2d}{d+1}}\frac{\lambda}{2}\right]}\leq \ivol[k]{Q^r} .
\end{equation}

Finally, notice that (\ref{part 1}), (\ref{useful condition}), and (\ref{part 2}) lead to the desired inequality $\ivol[k]{P^r}\leq \ivol[k]{Q^r}$ in a straightforward way. This completes the proof of Theorem~\ref{KP-type improved}.

\section{Proof of Theorem~\ref{improved contractions of molecules}}\label{second proof}

Clearly,
\begin{equation}\label{isodiametric-1}
{\rm diam}(Q_r)={\rm diam}(Q)+2r\leq\lambda+2r,
\end{equation}
where ${\rm diam}(\cdot)$ refers to the {\it diameter} of the corresponding set. Thus, the isodiametric inequality (Theorem 11.2.1 in \cite{BuZa}) and (\ref{isodiametric-1}) imply that
\begin{equation}\label{isodiametric-2}
\ivol[d]{Q_r}\leq\left(r+\frac{\lambda}{2}\right)^d\omega_d.
\end{equation}
Next, observe that (\ref{step 1 of defining f}), (\ref{step 2 of defining f}), (\ref{step 3 of defining f}), and (\ref{step 4 of defining f}) imply in a straightforward way that
\begin{equation}\label{density bounds}
\ivol[d]{P_{\mu\frac{\lambda}{2}}}\geq \frac{1}{f(d,\mu)}N\left(\frac{\lambda}{2}\right)^d\omega_d.
\end{equation}
For the next estimate recall that the {\it volumetric radius} of the compact set $\emptyset\neq A\subset \Ee^d$ is denoted by $R_{\Ee^d}(A)$ and it is defined by $\left(R_{\Ee^d}(A)\right)^d\omega_d:=\ivol[d]{A}$. Using this concept one can derive the following inequality from the Brunn-Minkowski inequality in a rather straightforward way (Theorem 9.1.1 in \cite{BuZa}):
\begin{equation}\label{Brunn-Minkowski applied 1}
R_{\Ee^d}(A_{\epsilon})\geq R_{\Ee^d}(A)+\epsilon,
\end{equation}
which holds for any $\epsilon>0$. Now, observe that as $r>\mu\frac{\lambda}{2}$ therefore 
\begin{equation}\label{set identity}
P_r=\left(P_{\mu\frac{\lambda}{2}}\right)_{r-\mu\frac{\lambda}{2}}.
\end{equation} 
Combining (\ref{set identity}) with (\ref{density bounds}) and (\ref{Brunn-Minkowski applied 1}) yields
\begin{equation}\label{Brunn-Minkowski applied 2}
\ivol[d]{P_r}\geq \left(\left(\frac{N}{f(d,\mu)}\right)^{\frac{1}{d}}\frac{\lambda}{2}+r-\mu\frac{\lambda}{2}\right)^d\omega_d
=\left(r+\left(\left(\frac{N}{f(d,\mu)}\right)^{\frac{1}{d}}-\mu\right)\frac{\lambda}{2}\right)^d\omega_d
\end{equation}
Finally, notice that $N\geq f(d,\mu)(1+\mu)^d$ implies $\left(\frac{N}{f(d,\mu)}\right)^{\frac{1}{d}}-\mu\geq 1$ which together with (\ref{isodiametric-2}) and (\ref{Brunn-Minkowski applied 2}) complete the proof of Theorem~\ref{improved contractions of molecules}.

\bigskip

%\normalsize

\noindent K\'aroly Bezdek \\
\small{Department of Mathematics and Statistics, University of Calgary, Canada}\\
\small{Department of Mathematics, University of Pannonia, Veszpr\'em, Hungary\\
\small{E-mail: \texttt{bezdek@math.ucalgary.ca}}

\end{document}